\documentclass[11pt]{amsart}
\usepackage{enumerate,amsmath,bm,amssymb,fleqn}
\newtheorem{discussion}{\bf Discussion}
\newtheorem{remark}{\bf Remark}
\newtheorem{proposition}{\bf Proposition}
\newtheorem{theorem}{\bf Theorem}

\newtheorem{example}{Example}
\begin{document}
\title[Study of spectral properties of dilation operators]{Effect of change in inner product on spectral Analysis of certain linear operators in Hilbert spaces arising from orthogonal polynomials}

\author{Anu Saxena}
\address{Department of Mathematics, Jesus and Mary College (University of Delhi), Delhi 110021, India. }
\email{asaxena@jmc.du.ac.in}
\maketitle
\begin{abstract}
In this paper we present how spectral properties of certain linear operators vary when operators are considered in different Hilbert spaces having common dense domain as the space of polynomials in one real variable with complex coefficients. This is done taking differential operator and matrix operator representation of dilation operator $E_{p,d}$ dilating $p$, a polynomial sequence, by $d$, a  non-constant sequence of non-zero complex scalars. For the purpose of identification of $E_{p,d}$ with formal differential operator, which includes finite order differential operators as special cases, we derive conditions under which a formal differential operator has a polynomial sequence as a sequence of eigen functions corresponding to non-zero eigen values. As a consequence we get, for instance, 

$$ ((1/72)x^4-3x) y^{(4)}-x y^{(2)}-(1/3 )x y^{(1)} + (4/3)y.$$

does not have any polynomial sequence as such a sequence of eigen functions, since no fourth order polynomial is such an eigen svector.\\
\textbf{Keywords:} Dilation Operator, Formal Linear Differential Operator, Matrix Operator, Ordinary Differential Operator, Shift Operator, Orthogonal Polynomials, Closure of an Operator, Spectrum, Eigen Functions.\\
\textbf{2010 AMS Subject classification:} 33C45, 34L05, 34L10, 47A05, 47A10, 47A58, 47A75
\end{abstract}
\section{Introduction}
Let $\mathcal{P}_c$ be the linear space of all polynomials in one real variable with complex co-efficients ; $\mathcal{P}$ be the set of sequences $p=(p_n)_{n\in \mathbb{N}_0}$, $p_n\in \mathcal{P}_c$, $\forall$  $n\in \mathbb{N}_0 $,  $p_0 \equiv 1$ and $p_{n}(x)=\sum_{k=0}^{n}p_{nk}x^k$,  $p_{nn} \ne 0$,  $n \ge 1$ ;  $\mathbb{N}_0$ is the set of non-negative integers. Let $\mathcal{P}_r$ be the real linear space consisting of real polynomials in one real variable. We take $\mathcal{P}_o$ = $\{ p \in \mathcal{P}: p \ \text {is an orthogonal polynomial sequence}\}$.
Any $p$ $\in \mathcal{P}$ is called a PS and  OPS stands for orthogonal polynomial sequence ; orthogonality is always with respect to a positive definite moment functional as in [1]. We denote by $\widetilde D$ the set of non-constant sequences $d=(d_n)_{n\in \mathbb{N}_0}$ of non-zero complex scalars. Set  $ D= \{ d \in \widetilde{D} :$ d is a real sequence with $ d_0=1\}$. \\
In this paper we consider dilation operators, in particular $S_{p,d}$ as defined and studied in [3]. Definition of $S_{p,d}$ was motivated by similar operators in literature , for instance, in [5] and [6]. For an OPS, $p \in\mathcal{P}_o $ and $d \in D$, $S_{p,d}$  is the linear map on $\mathcal{P}_c$ to itself such that  $S_{p,d}(p_n)= d_np_n$, $n \in \mathbb{N}_0$. Dilation map $E_{p,d}$ is like $S_{p,d}$  with  $p$ as a PS and $d \in \widetilde{D}$. For more facts on $S_{p,d}$ and $E_{p,d}$, we refer the reader to [3], [7] and to Appendix C. \\
For convenience of readers, relevant concepts and results needed throughout the paper have been freely collected from standard sources such as [1], [4], [8], [9] and research papers such as [3], [5], [6] and [7] and then, put as Appendices A to D given at the end of the paper.\\
We begin with giving some results from [5] and [6] in section 1. These motivated us to study and consider different representations of dilation operators.\\
In sections 2 and 3 we discuss and derive conditions when a formal differential operator, denoted as $\eta$ in this paper, has a dilation operator representation. As is well-known, for instance refer to [5] and Appendix B.1(i),
 $\eta$ : $\mathcal{P}_c \longrightarrow \mathcal{P}_c$ is the operator defined by the formal sum $$ \eta{(y)}(x) = \sum_{k=0}^{\infty} M_k(x) y^{(k)}(x)$$ where $ M_k(x)\in \mathcal{P}_c$ is such that $\deg M_k(x)\le k$ for
$k\in \mathbb{N} $, $M_0(x)$ is a constant, $y^{(0)}=y(x)$, $y^{(k)}(x)$ is the $k^{\rm th}$ derivative of $y(x)\in \mathcal{P}_c$ for
$k\ge 1$. Let $m : \mathcal{P}_c \longrightarrow \mathcal{P}_c $ be the operator $\eta$ for which $ M_{k}(x)\equiv 0$, $\forall$ $k>r$ for some smallest number $r \in \mathbb{N}$ with $M_r \not\equiv 0$. The number $r$ is then referred to as the order of finite order differential operator $m$. $\widehat{\eta}$ is an operator $\eta$ having a PS as a sequence of eigen functions corresponding to a sequence of non-zero eigenvalues. Similar meaning is attached to $\widehat{m}$.\\
In section 4 differential operator representation of dilation operator $S_{p,d}$ enables us to obtain a result about shift operator $\tau_{a,b}, a,b\in\mathbb{R}, a\neq 0$, defined as $\tau_{a,b}(x^n)=(ax+b)^n, n \in \mathbb{N}_0$, extended linearly to $\mathcal{P}_c$. It turns out that if such an operator has a sequence of orthogonal polynomials as a sequence of eigen functions then these polynomials are necessarily translates of some symmetric polynomial sequence with the corresponding sequence of eigen values just as $((-1)^n)$.\\
In section 5 dilation operator representation of differential operator enables us to do spectral analysis, in certain cases, by applying the Thin infinite matrix theory developed by the author and Ajit Iqbal Singh in [7]. Appendix D contains relevant definitions and results from [7]. In this section the setting changes from $\mathcal{P}_c$ to the Hilbert Space $H(q)$ arising from a PS $q$ in the following manner.
Starting with a PS $q=(q_n)$, we identify a complex sequence $a=(a_{n})_{n\in \mathbb{ N}_0}$, having finitely many non-zero terms, with $\sum_{n\in \mathbb{N}_0}a_n q_n$. Then completion of $\mathcal{P}_c$ with the inner product induced by the usual inner product in $l_2$ is a Hilbert Space which we denote by $H(q)$. If $q$ is taken as an OPS, then $\tilde{q}$ is the corresponding orthonormal sequence. We may note these spaces include Lebesgue Spaces $\mathcal{L}_\omega^{(2)}(a,b)$ whenever $q$ is an OPS with respect to weight function $\omega$ [1]. We begin with the question of closability of operators considered in this section. We then look at the spectral properties of operators which are closable and observe they display varied situations regarding spectra, closures and adjoints. In one instance $E_{p,d}$ in $H(q)$ has the whole of complex plane $\mathbb{C}$ in its spectrum, major part of which constitutes the continuous spectrum. There are instances where the domain of adjoint may be dense, yet no basic vector may be in it. Operators considered include finite order differential operators associated with classical orthogonal polynomials such as the ones given in Appendix B.3. We also observe that availability of closures, in some instances, is heavily dependent on eigen sequence $d$ and in some others nature of $d$ does not matter much.\\
\section{Motivation and Discussion}
We begin with quoting two results, Proposition MSZ from [6), Lemma~KS from [5],
in both retaining the original notation as far as possible.\\
\medskip\noindent\textbf{Proposition MSZ ([6]).} {\it
For a degree preserving linear operator $S$, the following statements are equivalent.
\begin{enumerate}[\rm(i)]\itemsep5pt
\item $s=S(1)$ and $S(y)=sy\circ \tau_{a,b}$, $y$ belonging to $\mathcal{P}_r$, $s\neq 0$,
$a\neq 0$, $b\in \mathbb{R}$.
\item If $(p_n)^\infty_{n=0}$ satisfies the three-term recurrence relation
\begin{align*}
xp_n(x)=a_np_{n+1}(x)+b_np_n(x)+c_np_{n-1}(x),
\end{align*}
then $(Sp_n)_{n=0}^\infty$ satisfies the three-term recurrence relation with \\
co-efficients $\alpha_n=a^{-1}a_n$, $\beta_n=a^{-1}(b_n-b)$ and $\gamma_n=a^{-1}c_n$, respectively.
\end{enumerate}}
Following lemma from [5] strongly indicates differential operator representation of the dilation operator.

\medskip \noindent\textbf{Lemma KS ([5]). }{\it
Consider the equation
\begin{eqnarray}
(L(y))(x)=\sum\limits_{k=1}^\infty M_k(x)y^{(k)}(x)=\lambda y(x),
\end{eqnarray}
where each $M_n(x)$ is a polynomial of degree $\le n$
\begin{eqnarray}
 M_n(x)=m_{n,0}+m_{n,1}x+\ldots+m_{n,n}x^n\,.
\end{eqnarray}
Let
\begin{eqnarray}
 \lambda_n=nm_{1,1}+n(n-1)m_{2,2}+\ldots+n!m_{n,n}
\end{eqnarray}
and assume that
\begin{eqnarray}
 \lambda_m\neq \lambda_n\quad (m\neq n).
\end{eqnarray}
Then
\begin{enumerate}[\rm(i)]\itemsep5pt
\item For $\lambda\neq \lambda_n$ $(n=0,1,\ldots)$ the only polynomial satisfying \emph{(1)}
is $y\equiv 0$. For $\lambda=\lambda_n$ there is a polynomial solution of \emph{(1)}, unique
to within an arbitrary constant multiplier, and it is of degree $n$.

\item Conversely, let $(P_n(x))$ be an
arbitrary PS and let $(\lambda_n)$ be an arbitrary sequence of constants with
$\lambda_n\neq 0$, $\lambda_m\neq \lambda_n$ $(m\neq n)$. Then there exists a unique sequence
$(M_n(x))$ of polynomials with degree $M_n\le n$, such that \emph{(1)} has $(y_n=P_n(x))$ as
its PS of solutions.
\end{enumerate}}

\begin{remark}\rm
For the purpose of identification of $E_{p,d}$ with $\eta $ through the
lemma, we first note the following facts:
\begin{enumerate}[(i)]\itemsep5pt\mathindent 2em
\item Though $\lambda_0$ is not defined in Lemma KS [5], \textsc{wlog} we take it as zero throughout. Since
$M_0p_0=d_0p_0$, $p_0\equiv 1$, $\forall$ PS   $p$ we get $M_0=d_0$, $\eta $ is $L+M_0I$, consequently in place of
$\lambda_n$ in the lemma we have $d_n-d_0$,
$n\in \mathbb{N}$.
\item  Part (i) then implies, for an $n\ge 1$, $\lambda_n=0\Longleftrightarrow d_n=d_0$ and for
$m,n\in \mathbb{N}_0$, $m\neq n$, $\lambda_m\neq \lambda_n\Longleftrightarrow d_m\neq d_n$.

\item We recall from [7], as in Appendix C.1(i), for a fixed $p\in \mathcal{P}_o$, $d\in D$, $S_{p,d}=S_{q,d}$ may be true
for infinite number of $q$. We shall, therefore, consider the triplet $(\eta ,p,d)$ for
differential operator representation of $E_{p,d}$. Also,
we treat the question of existence and that
of uniqueness separately for any one of $(\eta ,p,d)$ given
$p\in \mathcal{P}$, $d\in\widetilde{D}$.

\item
\begin{enumerate}[(a)]\itemsep5pt
\item Given $(p,d)$, unique $\eta $ may be obtained recursively as
\begin{align*}
M_1p_1^{(1)}&=(d_1-d_0)p_1 , \\
M_kp_k^{(k)}&=-\sum\limits_{j=1}^{k-1}M_jp_k^{(j)}+(d_k-d_0)p_k,\quad k>1,\quad k \in \mathbb{N}.
\end{align*}
Expanded form of the equations are given in Remark 3.
\item From Appendix B.3.(iii), we see that $m^\alpha(L_n^\alpha)=(-2n+1)L_n^\alpha$, $n\in \mathbb{N}_0$.
Uniqueness of $\eta$ then implies
\[S_{L_n^\alpha,d}\equiv m^\alpha \]
on $\mathcal{P}_c$ for $d=(-2n+1).$ \\
Same holds for $m^{\alpha,\beta}$ and $m^h$, given in Appendix B.3(i),(ii).
\end{enumerate}
\item It follows from Lemma KS that a necessary condition on $d\in \widetilde{D}$ for $(\eta ,p,d)$ to be such that $\eta  p=dp$ is that
\begin{align*}
d_n-d_0=\sum\limits_{r=1}^n m_{rr}p(n,r),\quad n\in \mathbb{N}.
\end{align*}
\item
Above condition on $d$ will be assumed without stating it when a \\
solution to $\eta  p=dp$ is sought. Also, $\sum\limits_{r=1}^n m_{rr}p(n,r)\neq -d_0$ is
implied.
\end{enumerate}
\end{remark}

\section{Differential operator representation of $\bm{(p,d)}$-dilation $\bm{E_{p,d}}$}

In this section, we discuss the formal differential operator $\eta$ and obtain conditions under which such an operator has a PS as a sequence of eigen functions corresponding to a sequence of non-zero eigen values. We refer to [4], book  by A.M. Krall, instead of original sources ; chapter XVII of the book contains a survey of results related to formal infinite order differential operators and associated classical-type orthogonal polynomials. In Appendix B.2, we give $\eta$ having Generalized Laguerre type polynomials (defined and studied by T.H.Koornwinder) as a sequence of eigen functions. Following Proposition identifies all $\eta $ which have a PS as such a sequence of eigen-functions. This, in turn, identifies all differential operators $\eta$ which do not have a PS as a sequence of eigen functions corresponding to a sequence of non zero eigen values. We also follow it up by a discussion of $\eta$ and $m$.
To ensure existence of solutions to $\eta  p=dp$, it is enough to obtain solutions to
$\eta  p_n=d_n p_n$ through induction on $n$, since $p_n^{(k)}=0$, for $k>n$.
\textsc{wlog} we assume $p_{nn}=1$ in Proposition 1.
\begin{proposition}
Assume $(p_r)_{r=0}^{n-1}$ exist such that $\eta  p_r=d_r p_r$.
For $k\in \mathbb{N}$, let $M_k$, $R_{k-1}\in \mathcal{P}_c$ be such that $\deg M_k\le k$, $\deg R_{k-1}\le k-1$
and $M_k=m_{kk}x^k+R_{k-1}$. Also, let
$\sum\limits_{k=1}^n p(n,k)R_{k-1}x^{n-k}=\sum\limits_{j=0}^{n-1}\alpha _j^{n-1}p_j$,
 where $(\alpha_j^{n-1})_{j=0}^{n-1}$ is an $n$-tuple of
scalars. Then $\exists \ p_n=x^n+q_{n-1}$ for some polynomial $q_{n-1}$
with $\deg q_{n-1}\le n-1:\eta p_n=d_np_n$
$$
\Longleftrightarrow
$$
$\exists $ an $n$-tuple $(\beta_j^{n-1})_{j=0}^{n-1}$ of scalars such that
$ (d_n-d_j)\beta_j^{n-1}=\alpha_j^{n-1}$, $0\le j\le n-1$.
\end{proposition}

\proof
$\sum\limits_{k=0}^n M_kp_n^{(k)}=d_np_n$.

\vskip.5em\noindent
$\Longleftrightarrow$ $\sum\limits_{k=1}^n (m_{kk}x^k +R_{k-1})
\left({p(n,k)x^{n-k}}+q_{n-1}^{(k)}\right)=(d_n-d_0)(x^n +q_{n-1}).$

\vskip.5em\noindent
$\Longleftrightarrow$ (i) $d_n-d_0=\sum\limits_{k=1}^n p(n,k)m_{kk}$,

\vskip.5em\noindent
\ (ii) $\sum\limits_{k=1}^n p(n,k)R_{k-1}x^{n-k}=(d_n-d_0)q_{n-1}-\sum\limits_{k=1}^{n-1}M_kq_{n-1}^{(k)}$\,.

\vskip.5em\noindent$\Longleftrightarrow \exists$ an $n$-tuple $(\beta_j^{n-1})_{j=0}^{n-1}$ of
scalars such that $q_{n-1}=\sum\limits_{j=0}^{n-1}\beta_j^{n-1}p_j$ satisfies
$\sum\limits_{j=0}^{n-1}\alpha_j^{n-1}p_j=(d_n-d_0)\sum\limits_{j=0}^{n-1}\beta_j^{n-1}p_j-\sum\limits_{j=0}^{n-1}(d_j-d_0)\beta_j^{n-1}p_j$\,.

\vskip.5em\noindent$\Longleftrightarrow$ $(d_n-d_j)\beta_j^{n-1}=\alpha_j^{n-1},\qquad 0\le j\le n-1\,.$
\endproof

\begin{remark}{\rm
From above, we have
\begin{enumerate}[(i)]\itemsep5pt
\item $d_n-d_j=0=\alpha_j^{n-1}$, for some $0\le j<n$ implies non-uniqueness of solutions.
In particular, we have the following facts.
\begin{enumerate}[(a)]
\item If $d_n=d_0$ for some $n\ge 1$ and $p_n$ is a solution of $\eta  p_n=d_np_n$ then $p_n+c$
is a solution for arbitrary scalar $c$.
\item If $d_t=d_n$ for some $t\neq n$, \textsc{wlog}, $n>t$, $t,n\ge 1$
and $p_n,p_t$ are such that $\eta  p_n=d_np_n$,
$\eta  p_t=d_t p_t$, then $\eta  (p_n-kp_t)=d_n(p_n-kp_t)$, $k\in \mathbb{C}$.
\item If $d_n=(-1)^n, n \in \mathbb{N}_0$,
then all symmetric PS are solutions of $\eta  p=dp$. In particular, it is
so for $p=(P_n^{\alpha,\alpha})$, for any $\alpha>-1$; refer to Appendix A.1.(ii).
\end{enumerate}
\item $d_n-d_j=0$, $\alpha_j^{n-1}\neq 0$, for some $n\in \mathbb{N}_0$ and $j\in \mathbb{N}_0$ with
$j<n$, implies non-existence of a required type of solution and this is the only way by
which a given $\eta $ fails to have a PS as a sequence of eigen functions.
\end{enumerate}}
\end{remark}

\begin{remark}\rm
\begin{enumerate}[(i)]
\item Alternative to Remark 2(ii) above, we can obtain all $\eta $ (in particular ${m}$) which do not have a PS as a sequence of eigen-functions from the following recursion equations.

\item \textbf{We give an expanded form of recursion equations }$\eta p_n=d_np_n$, $n\in \mathbb{N}_0$.\\
We first note that for compatibility of equations (a)-(e) below it is essential that
\begin{center}$(d_n-d_0)\neq -M_0$.\end{center}
\vskip.5em
\centerline{For $n\in \mathbb{N}_0$, we have $\eta p_n=d_np_n$}

\centerline{$\Longleftrightarrow$}
\begin{enumerate}[(a)]
\item $(d_n-d_0)p_{nn}=m_{nn}p(n,n)p_{nn}+\sum\limits_{r=1}^{(n-1)}m_{rr}p_{nn}p(n,r)$,

\item $(d_n-d_0)p_{n\,(n-1)}=m_{n\,(n-1)}p(n,n)\,p_{nn}$

\mbox{}\qquad\qquad\quad\qquad\qquad$+\sum\limits_{r=1}^{n-1}[m_{rr}\,p_{n\,(n-1)} p(n-1,r)+\,m_{r\,(r-1)}p_{nn}p(n,r)]$,
\item for $2\le r\le n-2$

\vskip.5em
$(d_n-d_0)p_{nr}=m_{nr}\,p_{nn}p(n,n)+\sum\limits_{k=r}^{n-1}\left[ \sum\limits_{t=0}^{\min\{n-k,r\}}m_{k\,(r-t)}p_{n\,(k+t)} p(k+t,k)\right]$

\vskip.5em
\mbox{}\qquad\qquad\quad\qquad\qquad\qquad\qquad$ + \sum\limits_{s=1}^{r-1}\left[ \sum\limits_{t=0}^{\min\{n-r,s\}}m_{s\,(s-t)}p_{n\,(r+t)} p(r+t,s)\right]$,

\item $ (d_n-d_0)p_{n1}=m_{n1}p_{nn}p(n,n)$

\vskip.5em
\mbox{}\qquad\quad\quad\qquad\qquad$+\sum\limits_{r=1}^{(n-1)} [m_{r1}p_{nr}p(r,r)+m_{r0}p_{n\,(r+1)}p(r+1,r+1)]$,
\item $ (d_n-d_0)p_{n0}=m_{n0}\,p_{nn}\,p{(n,n)}+\sum\limits_{r=1}^{n-1}m_{r0}p_{nr}p(r,r)$.
\end{enumerate}
\end{enumerate}
\end{remark}

\begin{remark}\rm
As a simple application, we get
\begin{align*}
{m}(y)=\left(12x^4-3x\right)y^{(4)}-xy^{(2)}-\dfrac{1}{3}xy^{(1)}+\dfrac{4}{3}
\end{align*}
as one differential operator not having a PS as a sequence of eigen functions. No
fourth order polynomial is an eigen function with a corresponding non-zero eigen value.
\end{remark}

\begin{discussion}{\rm
Operator $\widehat{m}$, a finite order differential operator has been studied by many authors. In [2], a survey of results concerning classification of $\widehat{m}$ has been done. Finite order and infinite order differential operators are intrinsically very different in many respects. For one, sequence of non-zero eigen-values of a finite order differential operator, when it exists, is an unbounded one. Therefore, for  $d$ bounded, $E_{p,d}$ is represented by an infinite order differential operator. The following Proposition points
how to construct many such operators.}
\end{discussion}
We point out that Operator in Appendix B.2 is an infinite order differential operator and its sequence of eigen values is an unbounded one.
\begin{proposition}
Any finite number of perturbations in the sequence of non-zero eigen-values of a finite order differential
operator, with the corresponding sequence of eigen-functions as fixed, makes it necessarily of infinite order.
\end{proposition}
\proof Let
$(\eta ,p,d)$ with $\eta$ as some finite order differential operator be given: $\eta  p=dp$.
Let $d'\in \widetilde{D}$ be such that for some $n\in \mathbb{N}_0$,
$d_k=d'_k$, $0\le k<n$ and $d_n\neq d'_n$. Also assume that $d_k\neq d'_k$ only for some finite number of indices $k>n$.
Then unique $\widetilde{\eta }\approx (\widetilde{M}_k)_{k\in \mathbb{N}_0}$ corresponding to $(p,d')$ is such that for $0\le j<n$, $\widetilde{m}_{jj}=m_{jj}$,
$\widetilde{m}_{n,n}-m_{n,n}=(d'_n-d_n)/n!$ and for $t>0$,
\begin{align*}
\widetilde{m}_{n+t,n+t}-m_{n+t,n+t}=\dfrac{d'_{n+t}-d_{n+t}}{(n+t)!}-\sum\limits_{r=n}^{n+t-1}(\widetilde{m}_{r,r}-m_{r,r})\dfrac{1}{(n+t-r)!},
\end{align*}
where $\widetilde{M}_k(x)=\sum\limits_{t=0}^k\widetilde{m}_{kt}x^t$, $M_k(x)=\sum\limits_{t=0}^k m_{kt}x^t$.
Now, if $\widetilde{\eta }$ were also a finite order differential operator, then we would have
\begin{align*}
\widetilde{m}_{n+t,n+t}-m_{n+t,n+t}=0,\quad \forall \ t>s \quad \text{for some} \ \ s\in \mathbb{N}.
\end{align*}
Also, $d'_{n+t}-d_{n+t}=0$, $\forall \ t>u$ for some $u\in \mathbb{N}$, by assumption.

\noindent Therefore, $\sum\limits_{r=n}^{n+t-1}(\widetilde{m}_{r,r}-m_{r,r})\dfrac{1}{(n+t-r)!}$
should be
zero, $\forall \ \ t>\max \{s,u,1\}$ ; which is not possible.
\endproof

\begin{remark}\rm
From the above recursion equations we notice that
for an arbitrary $\eta$,
 just one change in $d$, say, at $k^{\text{ th}}$ stage leads to
change in each subsequent $M_k$.
\end{remark}

\section{Shift-operators and $\bm{S_{p,d}}$}

In this section we discuss how differential operator
representation of $S_{p,d}$ helps
in concluding facts about shift operators. Shift operator, $\tau_{a,b}$ for $(a,b)\neq (1,0)$,
is an infinite
order differential operator $\eta $ with $M_k(x)=[(a-1)x+b]^k/k!$, $k\in\mathbb{N}_0$.
We give below
a characterization of shift operators having an OPS as sequence of eigen vectors via
identification with $S_{p,d}$ and observe that it is a very restricted class.

\begin{theorem}Let $p\in\mathcal{P}_o$, $d\in D$, $a,b\in\mathbb{R}:a\neq 0$.
Then $S_{p,d}=\tau_{a,b}$  $\Longleftrightarrow$ $d_n=(-1)^n$, $n\in \mathbb{N}_0$,
$a=-1$, $b$-arbitrary and $p_n=q_n\circ \tau_{1,b/2}$, $n\in \mathbb{N}_0$, for some symmetric OPS $q$.
\end{theorem}

\proof $S_{p,d}=\tau_{a,b}$ on $\mathcal{P}_r$ for some $a,b\in \mathbb{R}$, $a\neq 0$

\vskip.5em
\noindent $\Longleftrightarrow$ $S_{p,d}(p_n)=\tau_{a,b}(p_n)$, $n\in \mathbb{N}_0$

\vskip.5em
\noindent $\Longleftrightarrow$ $d_np_{nn}=p_{nn}a^n$, $n\in \mathbb{N}_0$ and other corresponding equalities hold

\vskip.5em
\noindent $\Longleftrightarrow$ $d_n=a^n$, $n\in \mathbb{N}_0$ and other corresponding equalities hold.

\vskip.5em
Next, taking $S=S_{p,d}$, $s=1$ in Proposition MSZ, Section 2 we have
\begin{align*}
S_{p,d}=\tau_{a,b}
\end{align*}
$\Longrightarrow$
if $(a_n)$, $(b_n)$, $(c_n)$ are the recurrence co-efficients for $(p_n)_{n\in \mathbb{N}_0}$, then
 $(a_n/a)$, $((b_n-b)/a)$,
$(c_n/a)$ are the recurrence co-efficients for $(S_{p,d}(p_n))$

\vskip.5em
\noindent $\Longrightarrow$ $a_n \left(d_n-\dfrac{d_{n+1}}{a}\right)p_{n+1}(x)
+d_n\left(b_n-\dfrac{b_n}{a}+\dfrac{b}{a}\right)p_n(x)$

\vskip.5em
\qquad\hskip1.3in $+c_n
\left(d_n-\dfrac{d_{n-1}}{a}\right)p_{n-1} {(x)}=0$, $n\in \mathbb{N}_0$, $x\in\mathbb{R}$

\vskip.5em
\noindent$\Longrightarrow$ for $n\in \mathbb{N}_0$, $d_n-d_{n+1}/a=0$ as well as
$b_n(1-1/a)+b/a=0$ and for $n\in \mathbb{N}$,
$d_n-d_{n-1}/a=0$.

\noindent
$\Longrightarrow$  $d_n=a^n$ and for $a\neq 1, b_n=\dfrac{-b}{(a -1)}$, for $n\in \mathbb{N}_0$. For $n\in \mathbb{N}$,
$d_n=a^{-n}$.

\vskip.5em
\noindent $\Longrightarrow$ $a=-1$, $d_n=(-1)^n$, $b_n=b/2$, $n\in \mathbb{N}_0$.

\vskip.5em
Thus $S_{p,d}=\tau_{a,b}$ implies $d_n=(-1)^n$, $n\in \mathbb{N}_0$, $a=-1$ and\\ $xp_n(x)=a_np_{n+1}(x)+\dfrac{b}{2}p_n(x)+c_np_{n-1}(x)$,
$n\in \mathbb{N}_0$, $x\in \mathbb{R}$.

\vskip.5em
Now,

\vskip.5em
\noindent $xp_n(x)=a_np_{n+1}(x)+b/2\,p_n(x)+c_np_{n-1}(x)$,
$n\in \mathbb{N}_0$, $x\in \mathbb{R}$

\vskip.5em
\noindent $\Longrightarrow$ $(x-b/2)p_n(x)=a_np_{n+1}(x)+c_np_{n-1}(x)$,
 $n\in \mathbb{N}_0$, $x\in \mathbb{R}$

\vskip.5em

\noindent $\Longrightarrow$ $tp_n(t+b/2)$= $a_np_{n+1}(t+b/2)+c_np_{n-1}(t+b/2)$, $n \in \mathbb{N}_0$, $t \in \mathbb{R}$.

\vskip.5em
\noindent $\Longrightarrow$ $q_n=p_n\circ \tau_{1,b/2}$ is symmetric.

\vskip.5em
It is easy to see that the conditions are sufficient also.
\endproof

\section{Spectral properties of the differential operator $\widehat{\eta}$ and
some matrix operators identified with $E_{p,d}$}

Having identified all differential operators $\widehat{\eta }$ in Section 3 we investigate
spectral properties of such operators along with those of some matrix operators.
All known finite-order differential operators having an OPS as a sequence of eigen functions corresponding to a sequence of non-zero eigen values are examples of $\widehat{\eta }$. Throughout we treat $\widehat{\eta}$ as an operator in some $H(Q)$.

\begin{discussion}\rm
\begin{enumerate}[(i)]\itemsep4pt
\item  First, we note that one can obtain the adjoint $\widehat{\eta}^*$ for an
arbitrary $\widehat{\eta}\approx E_{p,d}$, $p\in \mathcal{P}$, $d\in\widetilde{D}$ completely in terms of the co-ordinates of an
element in $H(Q)$, $Q$ a PS, an OPS or an ONS and the sequence $d$. Also, closure
may be obtained in similar terms in some cases since, for a linear operator $S$ in $H(Q)$
having $Q_n$, $n\in \mathbb{N}_0$ in the domain of $S$,
\begin{align*}
g\in D (S^*)\Longleftrightarrow \sum\limits_{k=0}^\infty |\langle g, S(Q_k)\rangle |^2, \ \text{is finite.}
\end{align*}
For $g\in D(S^*)$, $S^*g=\sum\limits_{k=0}^\infty \langle g,S(Q_k) \rangle Q_k$.
\vskip.5em
\item As an application, we consider the following four classes of operators and observe that
each case presents a
radically different situation.
\begin{enumerate}[(a)]\itemsep4pt
\item $E_{L_{n}^\alpha,d}$ in $H(\widetilde{L}_n^{\alpha+1})$, $\alpha>0$
\item $E_{L_{n}^{\alpha+1},d}$ in $H(\widetilde{L}_n^{\alpha})$, $\alpha>-1$
\item $E_{L_{n}^\alpha,d}$ in $H(L_n^{\alpha+1})$, $\alpha>-1$
\item $E_{L_{n}^{\alpha+1},d}$ in $H(L_n^{\alpha})$, $\alpha>-1$.
\end{enumerate}
\item
For notational convenience, from here onwards in this section, we refer to any one of the
above operators as $T$. We will see in the theorems that follow that all the above
operators are closable with appropriate conditions on $d$; hence, $D(T^*)$ is dense. Yet, we observe that elements common
to the $D(T^*)$ and the space of polynomials may vary anywhere from zero to the full space. In case (a) all basic vectors are always present in the domain of $T^*$ independent of what $d$ is. Hence, in this case, we have the closure for all operators $E_{L^\alpha_{n},d}$ in$H(\widetilde{L}_n^{\alpha+1})$, in particular, for ${m}^\alpha$, $\alpha>0$ in
$H(\widetilde{L}_{n}^{\alpha+1})$. In contrast, in (b)~either all basic vectors are present or none is in the domain of $T^*$.
For matrix operator in (c), the situation is entirely different, only a few basic vectors may be in $D(T^*)$ or may be none. For operator in (d) the situation is same as that in (b).
\item It may be emphasized that the same analysis may be
applied to many more differential operators satisfying the basic requirement of explicit
availability of connection co-efficients.
\end{enumerate}
\end{discussion}

\begin{theorem} Let $T$ be $E_{L_{n}^\alpha,d}$ in $H(\widetilde{L}_n^{\alpha+1})$,
$d\in \widetilde{D}$, $q=(L_n^{\alpha+1})$, $\alpha>0$. Let $g\in H(\widetilde{q})$. Set
$\ell=\sum\limits_{k=0}^\infty(g_k/r_k^{\alpha+1})$ with $r_k^{\alpha +1}=\|L_k^{\alpha +1}\|$,
\text{usual norm}.\\
\text{Then we have the following}

\begin{enumerate}[\rm(i)] \itemsep5pt

\item $g\in D(T^*)\Longleftrightarrow \sum\limits_{k=0}^\infty \left|g_k\bar{d}_k+\sum\limits_{t=0}^{k-1}(\bar{d}_t-\bar{d}_{t+1})\dfrac{r_t^{\alpha+1}}{r_k^{\alpha+1}}g_t\right|^2$ is finite\,,

\item for $g\in D(T^*)$, $T^* g=\sum\limits_{k=0}^\infty \left(g_k\bar{d}_k+\sum\limits_{t=0}^{k-1}(\bar{d}_t-\bar{d}_{t+1})\dfrac{r_t^{\alpha+1}}{r_k^{\alpha+1}}g_t\right)\widetilde{q}_k$\,,

\item $\widetilde{q}_s\in D(T^{*})$, $s\in \mathbb{N}_0$\,,

\item $T$ is closable and

$g\in D(\overline{T})\Longleftrightarrow \sum\limits_{s=0}^\infty \left|g_sd_s+\sum\limits_{k=s+1}^\infty (g_k/r_k^{\alpha+1})r_s^{\alpha+1}(d_s-d_{s+1})\right|^2$ is convergent
$\Longleftrightarrow$

$\sum\limits_{s=0}^\infty \left|g_sd_s+r_s^{\alpha+1}(d_s-d_{s+1})\left(\ell-\sum\limits_{k=0}^s (g_k/r_k^{\alpha+1})\right)\right|^2$ is finite\,,

\item for $g\in D(\overline{T})$, $\overline{T}g=\sum\limits_{s=0}^\infty \left(g_sd_s+r_s^{\alpha+1}(d_s-d_{s+1})\left(\ell-\sum\limits_{k=0}^s(g_k/r_k^{\alpha+1})\right)\right)\widetilde{q}_s$\,,

\item $g\in D(\overline{m^\alpha})$ $\Longleftrightarrow$ $\sum\limits_{s=0}^\infty \left| g_s(-2s+1)+2r_s^{\alpha+1}\left(\ell-\sum\limits_{k=0}^s \dfrac{g_k}{r_k^{\alpha+1}}\right)\right|^2$
is finite. For $g\in D(\overline{m^\alpha})$, $\overline{m^\alpha}(g)=\sum\limits_{s=0}^\infty \left( g_s(-2s+1)+2 r_s^{\alpha+1}\left( \ell-\sum\limits_{k=0}^s \left( \dfrac{g_k}{r_k^{\alpha+1}}\right)\right)\right) \widetilde{q}_s$.
\end{enumerate}
\end{theorem}
\proof
From Matrix form in Example 1, Appendix C.3 , we have for $k\in \mathbb{N}_0$,
\[T(\widetilde{q}_k)=d_k\widetilde{q}_k+\sum\limits_{t=0}^{(k-1)}(d_t-d_{t+1})\dfrac{r_t^{\alpha+1}}{r_k^{\alpha+1}}\widetilde{q}_t.\]
This gives (i) and (ii).\\
Discussion 2 (i) , then gives
\vskip.5em\noindent
$\widetilde{q}_s\in D(T^*)\Longleftrightarrow \sum\limits_{k=0}^\infty\left|(\bar{d}_s-\bar{d}_{s+1})\dfrac{r_s^{\alpha+1}}{r_k^{\alpha+1}}\right|^2$ is finite. \\
Now $(1/r_k^{\alpha+1})_{k\in \mathbb{N}_0}\in \ell_2$, refer to Appendix A.3. \\
Therefore $\widetilde{q}_s\in D(T^*)$, $s\in \mathbb{N}_0$, making $D(T^*)$ dense.\\
 Hence $T$ is closable. Also,
$T^*(\widetilde{q}_s)=\bar{d}_s\widetilde{q}_s+\sum\limits_{k=s+1}^\infty (\bar{d}_s-\bar{d}_{s+1})\dfrac{r_s^{\alpha+1}}{r_k^{\alpha+1}}\widetilde{q}_k$.\\
Since for a closable operator $T,\overline{T}=T^{**}$
\text{(iv) and (v) follow.}\\
Part (vi) follows from Remark 1(iv)(b), since in this case, $d_s-d_{s+1}=2$, $s\in \mathbb{N}_0$.

\begin{theorem}
Let $T$ be $E_{L_{n}^{\alpha+1},d}$ in $H(\widetilde{q})$, $q=(L^\alpha_{n})_{n\in \mathbb{N}_0}$, $\alpha>-1$,
 $d\in \widetilde{D}$. Let $g\in H(\widetilde{q})$. Then
\begin{enumerate}[\rm(i)]\itemsep4pt
\item $g\in D(T^*)\Longleftrightarrow \sum\limits_{k=0}^\infty \left|\bar{d}_kg_k+\sum\limits_{t=0}^{(k-1)}(\bar{d}_k-\bar{d}_{k-1})\dfrac{r_t^\alpha g_t}{r_k^\alpha }\right|^2$ is finite,

\item for $g\in D(T^*), T^* g=\sum\limits_{k=0}^\infty \left(\bar{d}_kg_k+\sum\limits_{t=0}^{(k-1)}(\bar{d}_k-\bar{d}_{k-1})g_t\dfrac{r_t^\alpha}{r_k^\alpha }\right)\widetilde{q}_k$,

\item for $s\in \mathbb{N}_0$, $\widetilde{q}_s\in D(T^*)\Longleftrightarrow ((\bar{d}_k-\bar{d}_{k-1})/r_k^\alpha )_{k\in \mathbb{N}_0}\in\ell_2$,

\item for $d\in \widetilde{D}:((\bar{d}_k-\bar{d}_{k-1})/r_k^\alpha)_{k\in N_0}\in\ell_2$,
\begin{enumerate}[\rm(a)]\itemsep4pt
\item $T$ is closable,

\item $g\in D(\overline{T})\Longleftrightarrow \sum\limits_{s=0}^\infty \left|g_sd_s+\sum\limits_{k=s+1}^\infty (d_k-d_{k-1})\dfrac{r_s^\alpha }{r_k^\alpha }g_k\right|^2$ is finite,

\item for $g\in D(\overline{T})$, $\overline{T}g=\sum\limits_{s=0}^\infty \left(g_sd_s+\sum\limits_{k=s+1}^\infty (d_k-d_{k-1})\dfrac{r_s^\alpha}{r_k^\alpha}g_k\right)\widetilde{q}_s $,
\end{enumerate}
\item for $\alpha>1$,
\begin{enumerate}[\rm(a)]\itemsep4pt
\item $ g\in D(\overline{m^{\alpha+1}})\Longleftrightarrow \sum\limits_{s=0}^\infty \left| g_s(-2s+1)+\sum\limits_{k=s+1}^{\infty}(-2)\dfrac{r_s^\alpha}{r_k^\alpha} g_k\right|^2$
is finite,

\item for $g\in D(\overline{m^{\alpha+1}})$,
$\overline{m^{\alpha+1}}(g)=\sum\limits_{s=0}^\infty \left( g_s(-2s+1)+\sum\limits_{k=s+1}^\infty (-2)\dfrac{r_s^\alpha}{r_k^\alpha} g_k\right)\widetilde{q}_s$\,.
\end{enumerate}
\end{enumerate}
\end{theorem}

\proof
It is enough to note that for $k\in \mathbb{N}_0$,
\begin{align*}
T(\widetilde{q}_k)=d_k\widetilde{q}_k+\sum\limits_{t=0}^{(k-1)}(d_k-d_{k-1})\dfrac{r_t^\alpha }{r_k^\alpha }\widetilde{q}_t\,,
\end{align*}
and for $d$ as in (iv), $T$ is closable by part (iii). The remaining parts may be proved as earlier.

 For part (v) we note that
from Remark 1(iv)b, we have $m^{\alpha+1}\equiv T$ in $H(L_n^\alpha)$ for $d=(-2n+1)$.
For this $d$, $d_k-d_{k-1}=\begin{cases}-2&,\ k\in \mathbb{N}\\ 1&,\ k=0\,.\end{cases}$

\noindent Therefore $m^{\alpha+1}$ is closable in $H(L_n^\alpha)$ if $\alpha>1$.
\endproof

\begin{theorem} Let $T=E_{p,d}$ be the operator defined via $p=(L^\alpha_{n})_{n\in \mathbb{N}_0}$
in $H(q)$, $q=(L^{\alpha+1}_n)_{n\in \mathbb{N}_0}$, $\alpha>-1$,
$d\in \widetilde{D}$. Let $g\in H(q)$. Then
\begin{enumerate}[\rm(i)]\itemsep4pt
\item $g\in D(T^*)\Longleftrightarrow \sum\limits_{k=0}^\infty \left|g_k\bar{d}_k+\sum\limits_{t=0}^{(k-1)}(\bar{d}_t-\bar{d}_{t+1})g_t\right|^2$ is finite,

\item for $g\in D(T^*)$, $T^* g=\sum\limits_{k=0}^\infty \left[g_k\bar{d}_k+\sum\limits_{t=0}^{(k-1)}(\bar{d}_t-\bar{d}_{t+1})g_t\right]q_k$,

\item for a fixed $j\in \mathbb{N}_0$, $q_j\in D(T^*)\Longleftrightarrow \bar{d}_j-\bar{d}_{j+1}=0$,

\item $E_{L^\alpha_{n},d}$ is unbounded $\forall$ $d\in \widetilde{D}$ as an operator in $H(L_n^{\alpha+1})$.

\end{enumerate}
\end{theorem}

\proof
From $T(q_k)=d_kq_k+\sum\limits_{t=0}^{(k-1)}(d_t-d_{t+1})q_t$, we have the first three parts following.

\noindent (iv) By definition $d\in \widetilde{D}$ is assumed to be a non-constant sequence. Part (iii)
above, therefore, gives that atleast one $q_j,j\in \mathbb{N}_0$, does not belong to $D(T^*)$.
Therefore, $D(T^*)\varsubsetneq H(q)$, $\forall$ $d\in\widetilde{D}$. Now, if $T$ is bounded for some
$d\in\widetilde{D}$, then $T$ is closable with $D(\overline{T})$ as the whole space $H(q)$ and $\overline{T}$ bounded. This, in turn,
implies that $D(T^*)=D(\overline{T}^*)$ as the full space $H(q)$. But this is contradictory. Therefore, $T$ is unbounded for every
$d\in\widetilde{D}$.
\endproof
\begin{remark}\rm
Matrix operators in the above theorem are closable for many $d$'s ; for instance, when $d_{k}-d_{k+1}= d_{k+2}-d_{k+3}, k \in \mathbb{N}_0 $, then $E_{p,d}$ is closable being thin. In such cases the spectrum is the whole of complex plane with major part as the continuous spectrum. \end{remark}

\begin{theorem} Let
$T=E_{p,d}$, in $H(q)$, $p=(L_n^{\alpha+1})_{n\in \mathbb{N}_0}$,
$q=(L_n^{\alpha})_{n\in \mathbb{N}_0}$, $\alpha>-1$, $d\in \widetilde{D}$. Let $g\in H(q)$. Then
\begin{enumerate}[\rm(i)]\itemsep4pt
\item $g\in D(T^*)\Longleftrightarrow \sum\limits_{k=0}^\infty \left|g_k\bar{d}_k+(\bar{d}_k-\bar{d}_{k-1})\sum\limits_{t=0}^{(k-1)}g_t\right|^2$ is finite,

\item for $g\in D(T^*)$, $T^* g=\sum\limits_{k=0}^\infty \left(g_k\bar{d}_k+(\bar{d}_k-\bar{d}_{k-1})\sum\limits_{t=0}^{k-1}g_t\right)q_k$,

\item for $s\in \mathbb{N}_0$, $q_s\in D(T^*)\Longleftrightarrow (\bar{d}_k-\bar{d}_{k-1})_{k\in \mathbb{N}_0}$ is in $\ell_2$.,

\item for $d$ such that $(\bar{d}_k-\bar{d}_{k-1})_{k\in \mathbb{N}_0}\in\ell_2$,

\begin{enumerate}[\rm(a)]
\item
$g\in D(\overline{T})$
$\Longleftrightarrow$
$\sum\limits_{s=0}^\infty \left|g_sd_s+\sum\limits_{k=s+1}^\infty (d_k-d_{k-1})g_k\right|^2$ is finite,

\item for $g\in D(\overline{T})$, $\overline{T}g=\sum\limits_{s=0}^\infty\left(g_sd_s+\sum\limits_{k=s+1}^\infty (d_k-d_{k-1})g_k\right)q_s$.
\end{enumerate}
\end{enumerate}
\end{theorem}

\proof Parts (i)-(iv) may be proved as earlier. We may note that by Theorem 3.3 in [7], given as Appendix item D.3 and Matrix form in Example 2, Appendix C.3,  $T$ is closable.
\endproof

\begin{remark}\rm
\begin{enumerate}[(i)]\itemsep4pt
\item For a closable operator $T$, $\sigma_0(\overline{T})=\{\lambda\in
\mathbb{C}:E_\lambda^{(2)} \not\equiv o\}$, where $E_\lambda^{(2)}$ is
the $\lambda$-approximate-Hilbert eigenspace of $T$ as defined in [7] given as
Appendix item D.2.
Here we note that in case $\overline{T}$
is completely determinable as above,
$(\overline{T}-\lambda I)g=0\Longleftrightarrow t_s-\lambda
g_s=0,s\in \mathbb{N}_0$, where $\overline{T}g=\sum\limits_{k=0}^\infty
t_kq_k$. For instance, in Theorem 5, $\lambda\neq d_k$, $k\in \mathbb{N}_0$,
\begin{align*}
(\overline{T}\!-\!\lambda I)g=o\Longleftrightarrow g_s
=-\!\!\sum\limits_{k=s+1}^\infty\!\dfrac{(d_k\!-\!d_{k-1})}{(d_s\!-\!\lambda)}g_k, \
s\in \mathbb{N}_0\,.
\end{align*}
\item  It may be emphasized that every PS or OPS serves as a sequence of eigen functions for an infinite
number of formal infinite order differential operators. The adjoint operators for these $\widehat{\eta}$ may be
obtained in $\ell_2$ or in some weighted $\mathcal{L}^2$ spaces. Closures can, similarly be obtained subject to conditions on $d$ in many
cases.

\item As has been shown in Theorem 5 that when $d\in \widetilde{D}$ is
such that $(\bar{d}_k-\bar{d}_{k-1})_{k\in \mathbb{N}_0}$ belongs to
$\ell_2$, we have the closure fully determined.
In Theorem 6 and Theorem 7,we take $d\in\widetilde{D}$ such that $(d_k-d_{k-1})_{k\in\mathbb{N}_0}\notin
\ell_2$. Here we observe heavy dependence of closure of the operator on the sequence of eigen values $d$.
\end{enumerate}
\end{remark}

\begin{theorem} \textbf{A set of necessary conditions:} Let $p,q,T$ be as in Theorem 5, $d \in \widetilde{D}$.
If $(f,g)\in\overline{G(T)}$, then
\begin{enumerate}[\rm(a)]\itemsep4pt
\item $\exists$ \ $(h_{n,u})_{n,u\in \mathbb{N}_0}$ such that
\begin{enumerate}[\rm{(a}1{)}]\itemsep4pt
\item $h_{n,u}=0$, $\forall$ $u>n$,
\item for $u\in \mathbb{N}_0$, $\lim\limits_{n\to \infty} h_{n,u}=f_u$,

\item $\lim\limits_{n\to\infty}h_{n,n}d_n=0$,
\item $\lim_{n\to\infty}\sum\limits_{u=1}^n h_{n,u}(d_u-d_{u-1})=g_0-f_0d_0$,
\end{enumerate}
\item for $k\in \mathbb{N}$,
\begin{align*}
 g_k=g_0-f_0d_0+f_kd_k-\sum\limits_{u=1}^k f_u(d_u-d_{u-1})\,.
\end{align*}
\end{enumerate}
\end{theorem}
\proof
The statement is a special case of Theorem 3.2 [7], refer Appendix D.5. Proof is as done earlier and is based on (1), (2) and (3)
below
\begin{enumerate}[(1)]
\item For
\begin{align*}
h_n&=\sum\limits_{u=0}^n h_{n,u}q_u \ \text{in the linear span of } \ \ q_n\text{'s}\\
Th_{n}&=\sum\limits_{k=0}^n t_{n,k}q_k,
\end{align*}
where
\begin{align*}
t_{n,k}=\begin{cases}\sum\limits_{u=k+1}^n h_{n,u}(d_u-d_{u-1})+h_{n,k}d_k&, \ \ 0\le k\le (n-1),\\
h_{n,n}d_n&,\ \ k=n\,.\end{cases}
\end{align*}
\item $(f,g)\in \overline{G(T)}\Rightarrow \exists \ (h_n)$ in
$H(q)$ such that $(h_n)\to f$ and $(Th_n)\to g$ in $H(q)$.

\item $T$ is closable by the above Theorem.
Therefore, $g$ is unique.\endproof
\end{enumerate}

\begin{theorem} \textbf{A set of sufficient conditions :} Let $p,q,T$ be as in Theorem 5, $d \in \widetilde{D}$.
Let $f\in H(q)$ be such that
\begin{enumerate}[\rm(i)]\itemsep4pt
\item $\lim\limits_{n\to\infty}\left[\sum\limits_{u=1}^n
f_u(d_u-d_{u-1})\right]$ is finite and, say, equals $S$,

\item for $(g_k)_{k\in \mathbb{N}_0}$ as $g_0=S+f_0d_0$,
\begin{align*}
g_k=S-\sum\limits_{u=1}^k f_u(d_u-d_{u-1})+f_kd_k,\quad k\ge 1,
\end{align*}
$(g_k)_{k\in \mathbb{N}_0}$ is in $\ell_2$,

\item $\lim\limits_{n\to\infty}(n+1)|f_nd_n-g_n|^2=0$.
\end{enumerate}
Then $f\in D(\overline{T})$ and
$g=\lim\limits_{n\to\infty}\left(\sum\limits_{t=0}^n g_tq_t\right)$ is
the unique element of $H(q)$ such that $(f,g)\in G(\overline{T})$.
\end{theorem}
\proof
Define $r_{n,u}$, $h_{n,u}$ for $n,u\in \mathbb{N}_0$ as
\begin{align*}
r_{n,u}&=\begin{cases}\frac{1}{n^22^n(|d_u-d_{u-1}|+|d_u|+1)}&,\ \ 0\le u\le n\\
-f_u&,\ \ u>n\end{cases}\\
h_{n,u}&=f_u+r_{n,u}.
\end{align*}
Then
\begin{align*}
h_n=\lim\limits_{k\to\infty}\left(\sum\limits_{u=0}^k h_{n,u}q_u\right)=\sum\limits_{u=0}^nh_{n,u}q_u
\end{align*}
is in $\mathcal{P}_c$, $\forall \ n\in \mathbb{N}_0$. Also,
$\lim\limits_{n\to\infty}\|h_n-f\|^2=0$.

Let $n\in \mathbb{N}_0$. Express $T(h_n)$ as in part (1) in the proof of Theorem 6.
Now,
\begin{align*}
\sum\limits_{k=0}^n |t_{n,k}-g_k|^2
&\le \left|\sum\limits_{u=1}^nr_{n,u}(d_u-d_{u-1})\right|^2 +\left|\sum\limits_{u=n+1}^\infty f_u(d_u-d_{u-1})\right|^2+|r_{n,0} d_0|^2\\
&\quad +\sum\limits_{k=1}^{(n-1)} \left|\sum\limits_{u=k+1}^n r_{n,u} (d_u-d_{u-1})\right|^2 +\sum\limits_{k=1}^{(n-1)}\left|\sum\limits_{u=n+1}^\infty f_u(d_u-d_{u-1})\right|^2\\
&\quad +\sum\limits_{k=1}^{(n-1)}|r_{n,k}d_k|^2+|r_{n,n}d_n|^2+\left|\sum\limits_{u=n+1}^\infty f_u(d_u-d_{u-1})\right|^2\\
&\le(n+1)\left|\sum\limits_{k=n+1}^\infty f_k(d_k-d_{k-1})\right|^2\\
&\quad + \text{remainder terms from definition of $r_{n,k}$}.
\end{align*}
The remainder terms tend to zero as $n$ tends to infinity by the definition of
$r_{n,k}$. Hence, from the assumptions, we have
\begin{align*}
0
&\le\limsup_{n\to\infty}\sum\limits_{k=0}^n |t_{n,k}-g_k|^2\\
&\le\limsup_{n\to\infty}(n+1)\left|\sum\limits_{k=n+1}^\infty f_k(d_k-d_{k-1})\right|^2\\
&=\lim_{n\to\infty}(n+1)|f_nd_n-g_n|^2\\
&=0
\end{align*}
$\Rightarrow \lim\limits_{n\to\infty}\|Th_n-g\|^2=0$. Hence,
$\exists \ (h_n)_{n\in \mathbb{N}_0}\in \mathcal{P}_c: (h_n)\to f$, $(Th_n)\to g$
implying $(f,g)\in \overline{G(T})$. This $g\in H(q)$ is unique,
since $T$ is closable.
\endproof

\subsection*{Acknowledgments}
I thank Professor Ajit Iqbal Singh for continued support and insightful suggestions.

\appendix
\section{\\Notation and terminology } \label{App:AppendixA}
% the \\ insures the section title is centered below the phrase: AppendixA
\def\N{\mathb{N}}
\def\cl{\mathop{\rm cl}}
\begin{itemize}
\item Let $\mathbb{N} _0$ be the set of non-negative integers and $ \mathbb{N} $, the set of natural numbers. Let $ \mathbb{R}$ be the set of Real numbers and $\mathbb{C}$ be the set of complex numbers.
\item $p(n,r)= \begin{cases} \frac{n!}{(n-r)!}, & n \ne 0\\  0 & n=0 \\ \end{cases}.$
\item A function is usually denoted by $f$, $g$, $y$  etc. but at times, in keeping with notation in a source, we write them as $f(x)$, $g(x)$, $y(x)$ etc. The context makes it clear.

\item For $n \in  \mathbb{N}_0$, let $(e_n)$ be the sequence $(\delta_{nm})_{m\in  \mathbb{N}_0}$.
\item Let $w$ be the space of all complex sequences and $\varphi$ be the subspace of
sequences with only finitely many non-zero terms. We equip $\omega$ with the topology
of pointwise convergence which is also given by the translation invariant metric
\[
d(x,y) = \sum_{n\in \mathbb{N}_0}\frac{|x_n-y_n|}{2^n(1+|x_n-y_n|)}\,.
\]
Then $(\omega, d)$ is complete.

\item Let $l_2$ be the subspace of $\omega$ consisting of square-summable sequences which becomes a Hilbert space under the usual inner product $\langle x, y\rangle =\sum_{n\in \mathbb{N}_0}x_n \bar y_n$.
\item We denote by $\widetilde{D}$ the set of sequences $d= (d_n)_{n \in \mathbb{N}_0}$ of non-zero complex scalars, $d\in \widetilde{D}$
is assumed to be non-constant to avoid triviality. Let $D $= \{ $ d\in\widetilde{D}: d$  is  a  real sequence with $d_0=1 $\}. Set $d_{-1} =0$. For $d \in \widetilde{D}$, closure of the set $\{d_n : n \in \mathbb{N}_0 \}$ is denoted by $cl(d)$.
\end{itemize}
\subsection{Some specific orthogonal polynomial sequences}
\begin{enumerate}[(i)]
\item For $\alpha > -1$, $(L_n^{(\alpha)})_{n\in  \mathbb{N}_0}$:
the Generalised Laguerre polynomial sequence is given by $L_n^{(\alpha)}(x) =
\sum_{k=0}^n \dfrac{(-1)^k}{k!}\binom{n+\alpha}{n-k} x^k$, where $\binom{t}{0} = 1$
and $ \binom{t}{k}  =\dfrac{t(t-1)\ldots(t-k+1)}{k!}$  for $t , k \in  \mathbb{N}$.
\item For $-1<\alpha < \infty , -1 < \beta < \infty $, $(P_{n}^{\alpha,\beta})_{n=0}^{\infty}$ : the Jacobi Polynomial Sequence is given by $P_{n}^{\alpha, \beta}(x) =1/2^n \sum_{m=0}^n \binom{n+\alpha}{m} \binom {n+ \beta}{n-m}(x-1)^{n-m} (x+1)^{m}$.
\item $(T_n)_{n\in \mathbb{N}_0}$: the Tchebyshev polynomials of the first kind are given by
\[
T_n(x) = \cos(n\theta), \ \ x = \cos \theta,\quad 0 \le\theta \le \pi.
\]

\item $(U_n)_{n\in  \mathbb{N}_0}$: the Tchebyshev polynomials of the second kind are given by
$U_n(x) =\dfrac{\sin(n+1)\theta}{\sin\theta}$, $x=\cos \theta$, $0 <\theta  <\pi$
extended by continuity to $[-1, 1]$ so that
\[
U_n(-1) = (-1)^n(n+1)\quad\text{and}\quad U_n(1) = n+1.
\]
\end{enumerate}

\subsection{Connecting relations}

Given $p = (p_n)_{n\in  \mathbb{N}_0}$, a PS, we set $p_{-1} = 0$, for notational convenience.
\begin{enumerate}[(i)]
\item $L_n^{(\alpha)} = L_n^{(\alpha+1)}-L_{n-1}^{(\alpha)}$,  $n\in  \mathbb{N}_0$.
\item $L_n^{(\alpha+1)} =\sum_{j=0}^n L_j^{(\alpha)}$,  $n\in  \mathbb{N}_0$.

\item $2T_n = U_n - U_{n-2}$, $n \in  \mathbb{N}$.

$T_0 = U_0$.

\item $U_{2n} = T_0 + 2 \sum_{k=1}^n  T_{2k}$,  $n \in  \mathbb{N}$,

$U_0 = T_0$.

$U_{2n+1} = 2 \sum_{k=0}^n T_{2k+1}$, $n\in  \mathbb{N}_0$.
\end{enumerate}

\subsection{An application of Raabe's test to Laguerre polynomials}
 For Laguerre polynomials $L_k^\beta$, $\beta>-1$, $k\in\mathbb{N}_0$, we have the usual norm,
$r_k^\beta=\|L_k^\beta\|=\sqrt{\dfrac{\Gamma(k+\beta+1)}{p(k,k)\Gamma(\beta+1)}}$.
\smallskip
Thus
$r_k^\beta=\begin{cases}\sqrt{(1+\beta)(1+\beta/2),\ldots,(1+\beta/k)}&,\ k\in\mathbb{N}\\1&, \ k=0\,.\end{cases}$
If $\beta>1$, then by Raabe's Test $\left(\frac{1}{r_k^\beta}\right)\in\ell_2$.

\section{\\Formal Differential operator} \label{App:AppendixC}
% the \\ insures the section title is centered below the phrase: Appendix B
 \subsection {$\eta$ ; Formal differential operator. }

\begin{enumerate}[(i)]
 \item Let $\eta$ : $\mathcal{P}_c \longrightarrow \mathcal{P}_c$ be the operator defined by the formal sum $$ \eta{(y)}(x) = \sum_{k=0}^{\infty} M_k(x) y^{(k)}(x)$$ where $ M_k(x)\in \mathcal{P}_c$ is such that $\deg M_k(x)\le k$ for
$k\in \mathbb{N} $, $M_0(x)$ is a constant, $y^{(0)}=y$, $y^{(k)}$ is the $k^{\rm th}$ derivative of $y\in \mathcal{P}_c$ for
$k\ge 1$. For $p\in \mathcal{P}$, $d\in\widetilde{D}$ equations $\eta  p_n=d_np_n$, $n\in \mathbb{N}_0$ collectively
are abbreviated as $\eta  p=dp$.\\ For $n \in \mathbb{N}_0$, when $p_n$ is taken as a solution of $\eta  p_n=d_np_n$, $n\in \mathbb{N}_0$, it will be assumed that deg $p_n$=$n$. Also $d_n$ will be referred to as an eigen-value and $p_n$ as the corresponding eigen-function.
\item $m$ ; Let $m : \mathcal{P}_c \longrightarrow \mathcal{P}_c $ be the operator $\eta$ for which $ M_{k}(x)\equiv 0$, $\forall$ $k>r$ for some smallest $r \in \mathbb{N}$ with $M_r \not\equiv 0$. Number $r$ is then referred to as the order of finite order differential operator $m$.
\item $\widehat{\eta}, \widehat{m}$ ; $\widehat{\eta}$ is an operator $\eta$ having a PS as a sequence of eigen functions corresponding to a sequence of non-zero eigenvalues. Similar meaning is attached to $\widehat{m}$.\\
Study of arbitrary $m$ is implicit while studying $\eta$ in general. Similarly treatment of $\widehat{m}$ automatically comes under any analysis done for $\widehat{\eta}$.
\end{enumerate}
\subsection{ An infinite order differential operator}
We give below details of a well-known infinite order differential operator  $\eta$. We refer to [4], Chapter XVII for more information on the topic.

\item For $\alpha >-1,  K > 0, $ the Generalized Laguerre-type polynomials $(
L_{n}^{ \alpha,  K})_{n=0}^{\infty}$ are defined thus, $$ L_{n}^{\alpha, K}(x)= \left[1+K\binom{n+\alpha}{n-1}\right] L_{n}^{\alpha}(x)+ K\binom{n+ \alpha}{n}\frac{d}{dx}L_{n}^{\alpha}(x),n \in \mathbb{N}_0, $$  where $L_{n}^{\alpha}(x), n \in \mathbb{N}_0$ are the Generalized Laguerre polynomials.  Let $$ \eta_{\alpha, K}{(y)}(x) = \sum_{k=0}^{\infty} M_k(x) y^{(k)}(x)$$,
where
\begin{align*}
M_{k}(x)=\begin{cases}K \frac{1}{k!}\sum\limits_{j=1}^k (-1)^{k+j+1}\binom{\alpha + 1}{j-1} \binom{ \alpha +2}{k-j} (\alpha +3)_{k-j}x^k &, \ \ k \ge 2,\\
-Kx+ \alpha +1&,\ \ k=1, \\
1 &,   k=0. \\  \end{cases}
\end{align*}
Take $\alpha, K$ such that $d_n=-K\binom{n+ \alpha +1}{n-1}-n+1 \ne 0, n \in \mathbb{N}_0 .$ Let $p_n= L_{n}^{\alpha, K}, n\in \mathbb{N}_0.$ Then $\eta_{\alpha,K}(p_n)=d_np_n, n\in \mathbb{N}_0.$

\subsection{Differential Operators associated with some classical OPS }

\begin{enumerate}[(i)]
\item  Let $-1<\alpha < \infty, -1< \beta <\infty, \alpha + \beta \ne -1, p= (P_{n}^{\alpha, \beta})$. Let $m^{\alpha, \beta}$ be the differential operator determined by $$(m^{\alpha,\beta}f)(x)=(1-x^2)f^{\prime\prime}(x)+(\beta -\alpha-(\alpha+\beta+2)x)f'(x)+f(x).$$ For $$d_n=-n(n+\alpha+\beta+1)+1,   n\in \mathbb{N}_0$$  $$ m^{\alpha,\beta}(p_n) =d_np_n, n\in \mathbb{N}_0.$$ Condition $\alpha+\beta \ne -1$ is imposed to have $d_n \ne 0, n\in \mathbb{N}_0$.
\item Let $p=(H_n)_{n=0}^{\infty} $, the Hermite polynomial system,
$d=(-2n+1)_{n\in \mathbb{N}_0}$ and $m^h$ be the differential operator determined by $$(m^{h}f)(x)=f^{\prime\prime}(x)-2xf'(x)+f(x).$$ Then $m^{h}(p_n)=d_np_n, n \in \mathbb{N}_0$.
\item Let $\alpha >-1$ and $p= (L_{n}^{\alpha})_{n=0}^{\infty}$, the Generalised Laguerre system. Let $m^ {\alpha}$ be the differential operator determined by $$ (m^{ \alpha}f)(x)= 2xf^{\prime\prime}(x)+2( \alpha +1-x)f'(x)+f(x).$$ Then for $d=(-2n+1)_{n\in \mathbb{N}_0}$, $$m^{\alpha}(p_n)=d_np_n, n \in \mathbb{N}_0.$$
\end{enumerate}

\section{\\Hilbert Space $H(Q)$ and Dilation Operators} \label{App:AppendixC}

\subsection {($p,d$) dilation operators }
\begin{enumerate}[(i)]
\item  $S_{p,d}$ : For $d \in D$ and $ p \in \mathcal{P} $ let $S_{p,d}$ be the linear operator on $\mathcal{P}_c$ to itself given by $S_{p,d}(p_n)= d_np_n$, $n \in \mathbb{N}_0$ . For a fixed $p$,  $S_{p,d}$ and $d$ determine each other. On the other hand, for a fixed $d$, $p$ is not always determined by the operator $S_{p,d}$ as it may be equal to $S_{q,d}$ for an infinite number of $q\textsc{\char13}s$. One such $d$ is $d_n=(-1)^n, n \in \mathbb{N}_0$. For this $d$, $S_{p,d}$=$S_{q,d}$ , on $\mathcal{P}_c$ for all $\gamma > -1,  p= (P_{n}^{\gamma,\gamma } ),  q=(T_n) $.

\item  Dilation map $E_{p,d}$: When we relax the condition either on  $p$ or on $d$ in the definition of $S_{p,d}$ above to having $p$ as an arbitrary $p \in \mathcal{P}$, $d \in \widetilde{D}$, we denote the resulting operator by $E_{p,d}$\\
 \end{enumerate}
\subsection {Hilbert Spaces $H(q)$} Let $q = (q_n)$ be a PS. We consider $\mathcal{P}_{c}$ as the linear span
             of $\{e_n = q_n : n \in\mathbb{N}_0\}$ and identify it with $\varphi$ making
            $a = (a_n)_{n\in\mathbb{N}_0}\in\varphi$ correspond to $\sum_{n\in\mathbb{N}_0}a_nq_n$.
             Then the completion of $\mathcal{P}_{c}$ with the inner product induced by the usual
             inner product in $\varphi\subset l_2$ is a Hilbert space, which we denote
             by $H(q)$.
             \begin{enumerate}[(i)]
             \item If $q$ is taken as an OPS, then $\tilde{q}$ is the corresponding ONS. When referring to $q$ a PS, an OPS or an
             ONS, we shall write $\mathbf{ Q}$. \\An arbitrary $g\in H(Q)$  will be identified with the sequence $(g_j)_{j\in \mathbb{N}_0} \subset  l_2$ such that $g= \sum_{j\in \mathbb{N}_0}g_j Q_j$, where $Q=(Q_n)$. Elements of $Q=(Q_n)$ in $H(Q)$ will be referred to as basic vectors.\\
             \item When $q$ is an OPS orthogonal with respect to a weight function $w$ with the interval of orthogonality as $(a,b)$ and satisfying the normalizing condition $\int_{a}^b w(x)dx=1$, then $H(Q)$ is $\mathcal{L}_\omega^{(2)}(a,b)$. The interval of orthogonality may be infinite . For instance, for $\alpha > -1$, $ H({\widetilde L}_n^{\alpha})$ is $\mathcal{L}_\omega^{(2)}(0,\infty)$ with weight $w(x)= \frac {x^\alpha \exp {(-x)}}{\Gamma{(\alpha )}}$, for $x>0$.
            \item  We refer to [1] for orthogonal polynomials.
             \end{enumerate}

\subsection {Some matrix Operators}  The linear operator $T$ in $H(q)$ determined by $S_{p,d}$ on $\mathcal{P}$ is given by the matrix
$A = (a_{jk})$ in each of the following examples, when $p$ and $q$ are taken as specified.
\begin{example}
For $p =(p_n) = (L_n^\alpha)$ and  $q = (q_n)=(L_n^{\alpha+1})$, $\alpha > -1$, we have
$A = (a_{jk})$ with
\begin{align*}
&a_{jk}=0,\quad k<j,\\
&a_{jj}=d_j,\quad j\ge 0,\\
&a_{jk}=d_j-d_{j+1},\quad k>j\,.
\end{align*}
\end{example}
\begin{example}
For $p = (p_n) = (L_n^{\alpha+1})$ and $q = (q_n) = (L_n^\alpha)$, $\alpha > -1$, we have
$A = (a_{jk})$ with
\begin{align*}
&a_{jk} = 0,\quad k< j,\\
&a_{jj}=d_j,\quad  j \ge 0,\\
&a_{jk} = d_k- d_{k-1},\quad 0 \le j<k, \ k \ge 1 .
\end{align*}
\end{example}
\section{\\Thin Matrices and Closability of Associated Operators in Hilbert Spaces} \label{App:AppendixB}
% the \\ insures the section title is centered below the phrase: Appendix B

We give below relevant definitions, theorems and examples from the paper [7] by the author and Ajit Iqbal Singh. We shall confine our attention to the setting of Hilbert Space $l_2$ though the paper begins with that of Locally Convex Space $\omega$ before coming to the Hilbert Space setting.
An infinite complex matrix  $A=(a_{jk})_{j,k \in \mathbb{N}_0}$ determines an operator $T$ from $\varphi $ to $\omega$  in a natural way. Also, it is well-known that Columns of $A$ are in $l_2 \Longleftrightarrow $ $T$ takes $\varphi$ into $l_2$  .

\subsection{Definition of Thin infinite matrices}

\begin{enumerate}[(i)]
\item We define an equivalence relation $\sim$ in $\omega$ as : For $a$, $b$ in $\omega$,  $a\sim b$ if and only
    if $a - \mu b \in l_2$ for some non-zero $\mu$ in $\mathbb{C}$.

      The relation $\sim$ divides $\omega$ into mutually disjoint equivalence classes $[a]$
    with $[o] = l_2$. Further for $a\not\sim o$, $a\sim b$,  the constant $\mu$ is unique.

\item For an infinite matrix $A$ denote by $a_k$, the $k$th row, $k \in \mathbb{N}_0$. Then there exists an index set $I\subseteq \mathbb{N}_0$ and a decomposition $(N_i)_{i\in I}$ of $\mathbb{N}_0$ into mutually disjoint sets such that for $k, t \in \mathbb{N}_0$, $a_k \sim a_t$ if and only if $k, t \in N_i$ for some $i \in I$.

      For notational convenience, if $a_k \in l_2$ for some $k \in \mathbb{N}_0$ then we take
    N$_0 = \{t:a_k\sim a_t\}$ and otherwise we take $I\subset \mathbb{N} $. For $i \in I$, let
    $k_i^0 =\min\{k : k \in N_i\}$ and for $k \in N_i$ with $i \ge 1$, let $m_k$ be
    the unique non-zero number such that $a_k- m_ka_{k_i^0} \in l_2$.

 Clearly, $m_{k_i^0} =1$.

If $0 \in I$, then we set $m_t = 0$ for $t$ in $N_0$. Let $m = (m_k)_{k \in \mathbb{N}_0}$.

If for an $i \in I$, $i \ge 1$, we have $N_i$ as an infinite set, then we arrange $N_i$ as
a sequence, say, $(t_i)_{t\in \mathbb{N}_0}$ with $0_i = k_i^0$ and $t_i \neq t'_i$ for $t\neq t'$ in $\mathbb{N}_0$.

Next we put $b_{jk} = a_{jk}- m_ja_{k_i^0k}$ for $j \in N_i$, $i \in I$, $k\in \mathbb{N}_0$. We  say
       that $A$ is \emph{$(m, (k_i^0)_{i\in I})$-thin} with \emph{thinning} $B = (b_{jk}),{j,k\in  \mathbb{N}_0}$.

\item
\begin{enumerate}[(a)]
\item Clearly, $V = B$ restricted to $l_2$ is a continuous linear operator on $l_2$ into $\omega$.

\item If $I$ contains an $i\ge 1$, i.e., not all the rows of $A$ are in $l_2$, then $T$ is not
      continuous as an operator from $l_2$ to $\omega$. To see this well-known fact, we first note that $s_n =\sum_{k=0}^n |a_{k_i^0k}|^2\uparrow \infty$.

      Let $n_0$ be the least integer with $s_{n_0}\neq 0$. We next define the sequence
      $(h^{(n)})$ in $\varphi$ by taking $h_t^{(n)} =\dfrac{\overline{a_{k_i^0t}}}{s_n}$, $n \ge n_0$, $0 \le t\le n$,
      and $h_t^{(n)} = 0 $,
      otherwise. Then $h^{(n)}\to o$ in $l_2$ whereas, $(Th^{(n)})_{k_i^0} = 1$ for $n \ge n_0$.
\end{enumerate}
\item A will be called \emph{thin} if either $I = \{0\}$ or for each $i \in I$, $i \ge 1$, $N_i$ is infinite
       and the sequence $m^{(i)} = (m_{t_i})_{t\in \mathbb{N}_0}\notin l_2$.

\item A will be called \emph{blocked} if $a_{jk} = 0$ for $j$ and $k$ not in the same $N_i$ for any $i \in I$.

\end{enumerate}

\subsection{$\lambda$-approximate-Hilbert eigenspace}
Suppose columns of $A$ are in $l_2$. Then, for $\lambda \in \mathbb{C}$, $E_\lambda^{(2)} = \{x : (x,\lambda x) \in \overline{G(T)}\}$
will be called the \emph{$\lambda$-approximate-Hilbert eigenspace} of $T$.
\subsection{Criterion for closability} {Theorem 3.3 from [7]}:
Let $A$ be an infinite matrix whose columns are in $l_2$. Then the following hold.
\begin{enumerate}[\rm(i)]
\item If $A$ is thin then $T$ is closable.
\item If $A$ is blocked and $T$ is closable then $A$ is thin.
\end{enumerate}
\subsection{Examples of Thin infinite matrix, infinite blocked matrix}
\begin{enumerate}
\item Example.2 in Appendix C.3 is thin.
\item Matrix in this example is blocked.

Let $p = (p_n)$: \[
p_0 = T_0,\quad p_n = 2T_n\qquad\text{for} \ n \in\mathbb{N}_0, \ n \ge 1
\]
and $q = (q_n) = (U_n)$.

Then the matrix for $T=S_{p,d}$ in $H(q)$ is given by
$A = (a_{jk})$
\begin{align*}
&a_{jk} = 0,\qquad k < j\\
&a_{jj} = d_j,\quad j \ge 0\\
&a_{jk} = d_j-d_{j+2},\quad k=j+2s, \ s\in \mathbb{N}\\
&a_{jk} = 0,\qquad k = j + (2s - 1),\quad s\in\mathbb{N}.
\end{align*}
Therefore, $T=S_{p,d}$ can be thought of as $T'\oplus T''$ acting on $\mathcal{P}'\oplus \mathcal{P}''$  with
\begin{align*}
\mathcal{P}'&=\text{the linear span of $\{q_{2n}: n \in \mathbb{N}_0\}$ and}\\
\mathcal{P}''&=\text{the linear span of $\{q_{2n+1}: n \in \mathbb{N}_0\}$},
\end{align*}
$T'$ and $T''$ are given by matrices $A'$ and $A''$ respectively of the same form as
$A$ in Example 1 in Appendix C.3, with $d$ replaced by $d' = (d_{2n})_{n\in \mathbb{N}_0}$ and $d'' = (d_{2n+1})_{n\in \mathbb{N}_0}$.

Thus cardinality of $I \le 3$, where $I$ is an indexing set as in D.1(ii) above.
Hence Matrix $A$ is blocked.
\end{enumerate}

\subsection {Theorem 3.2 [7]}

Let $A$ be an infinite matrix whose columns are in $l_2$. Let $m$, $(k_i^0)_{i\in I}$
and $V$ be as in D.1.(ii) and (iii)  above.
\begin{enumerate}[\rm(i)]
\item Let $(x,y) \in \overline{G(T)}$. Then
\[
y_t = y_{k_i^0 } m_t+(Vx)_t,\quad t\in N_i, \ i\in I\,.
\]
\item If $o\neq x \in  E_0^{(2)}$ then $0$ is an eigenvalue of $V$ with $x$ as an eigenvector.

\item For $\lambda \neq 0$ in $\mathbb{C}$ and $o\neq x\in E_\lambda^{(2)}$, $x_{k_i^0}=0$ for all $i \in I$, $i\ge 1$ if and only if
      $\lambda$ is an eigenvalue of $V$ with $x$ as an eigenvector.
\item If $o\neq x \in E_\lambda^{(2)}$ then,
\[
(\lambda I- V)x = \lambda \big((x_{k_i^0} m_t)_{t\in N_i}\big)_{i\in I}
\]
and, thus, for $\lambda \neq 0$  in $\mathbb{C}$ we have $\big((x_{k_i^0} m_t)_{t\in N_i}\big)_{i\in I}\in R(\lambda I-V)$.

\item If $0$ is not an eigenvalue of $V$ then $T^{-1}$ exists and is closable.

\item If $\lambda (\neq 0)$ in $\mathbb{C}$ is not an eigenvalue of $V$ with an eigenvector in the
subspace $\{x: x_{k_i^0} = 0, i \in I, i \ge 1\}$ and for no non-zero sequence $(\alpha_{k_i^0})_{i\in I}$,
$((\alpha_{k_i^0}m_t)_{t\in N_i})_{i\in I}$ is in $R(\lambda I - V)$ then $(\lambda I- T)^{-1}$ exists and is closable.
\end{enumerate}

\end{document}